\documentstyle[12pt]{article}
\begin{document}

\centerline{\bf \LARGE On the $L^{n\over 2}$-norm of Scalar Curvature}

\vspace{0.3in}

\centerline{\bf \large LEUNG, Man Chun}
\vspace{0.075in}
\centerline{\bf \large Department of Mathematics,}
\centerline{\bf \large National University of Singapore,}
\centerline{\bf \large Lower Kent Ridge Road,}
\centerline{\bf \large Singapore 0511}
\centerline{\bf \large matlmc@leonis.nus.sg}

\vspace{0.35in}

\begin{abstract}
Comparisons on $L^{n\over 2}$-norms of scalar curvatures between
Riemannian metrics and
standard metrics are obtained. The metrics are restricted to conformal
classes or
under certain curvature conditions.
\\ \end{abstract}

\vspace{0.35in}

{\bf \LARGE 1. \ \ Introduction}

\vspace{0.3in}

Let $M$ be a compact $n$-manifold without boundary. For a Riemannian
metric $g$
on $M$, curvature tensor, Ricci curvature tensor and scalar
curvature of $g$ are denoted by $R (g)$, Ric $(g)$ and $S (g)$, respectively.
A natural and interesting question in Riemannian geometry is relations
between
topology of the manifold $M$ and curvatures of $g$. Often topology of $M$
would
impose certain restrictions on the behavior of curvatures of the metric
$g$.
The Gauss-Bonnet theorem provides a beautiful relation in this direction.
As complexity of the Gauss-Bonnet
integrand increases with dimension, it would be desirable to obtain
simpler but
not "sharp" relations. Indeed, there have been many
interests on $L^{n\over 2}$-curvature pinching
and bounds on topological quantities by integral norms of curvatures. In
this article, we study some questions on obtaining
lower bounds on $L^{n\over 2}$-norms of the Ricci curvature
and scalar curvature. There are some rather general and well-known
problems:
given a compact $n$-manifold $M$, for a sufficiently large class of
Riemannian metrics $g$ on $M$,
whether there are positive lower bounds on the following:
\\[0.25in]

(1) \ \ Vol $(M, g)\,,$ \ provided $K_g \ge -1$ \ or \ Ric $(g)_{ij} \ge
- (n - 1)
g_{ij}$ \ or \ $S (g) \ge - n ( n - 1)$, where $K_g$ is the sectional
curvature of $(M, g)$;\\[0.01in]
(2) \ \ $\int_M |S (g)|^{n \over 2} dv_g$;\\[0.01in]
(3) \ \ $\int_M | {\mbox {Ric}}\, (g) |^{n\over 2} dv_g\,.$\\[0.02in]
We note that (2) and (3) are both scale invariant, while a lower bound on
curvature is
required in (1) so as Vol $(M, g)$ would not go to zero by scaling. As a
flat torus would not have positive lower bounds on (1), (2) and (3),
some restrictions are needed on the manifold $M$.
Some suggestions are:\\[0.1in]
(a) \ \ $M$ admits a locally symmetric metric of strictly
negative sectional curvature, or\\[0.01in]
(b) \ \ $M$ admits an Einstein metric of negative sectional
curvature, or simply\\[0.01in]
(c) \ \  $M$ admits a metric of negative sectional curvature.\\[0.02in]
\hspace*{0.3in}Recently, Besson, Courtois and Gallot [5,6] have
demonstrated that if $(M, h)$ is
a compact hyperbolic $n$-manifold ($n \ge 3$), then for any Riemannian
metric $g$ on $M$ with
Ric $(g) \ge -(n - 1) g$, one has Vol $(M, g) \ge {\mbox {Vol}}\, (M, h)$
and equality holds if and only
if $(M, g)$ is isometric to $(M, h)$. In this note, we would mainly
consider question (2) and
(3), under one of the conditions in (a), (b) or (c) and with restrictions
on the
choices of the Riemannian metric $g$ by certain curvature assumptions or
in a certain
conformal classes. Our method is to investigate relations between
the $L^{n\over 2}$-norms
of scalar curvatures for different metrics with that of a standard
metric.\\[0.02in]
\hspace*{0.3in}The Gauss-Bonnet theorem for two-manifolds shows that, if
$M$ is a compact
surface and $h$ is a metric on $M$ with constant negative curvature $S
(h)$, then
$$\int_M |S (g)| dv_g \ge \int_M |S (h)| dv_h\,. \leqno (1.1)$$
Let $\chi (M)$ be the Euler characteristic  of $M$.
The Gauss-Bonnet theorem for higher dimensions ($n$ even) states that [18]
$$c_n \chi (M) = \int_M \sum_{\sigma \in {\cal C}_n} \sum_{\tau \in {\cal
C}_n}
\varepsilon (\sigma) \varepsilon (\tau)
 R (g)_{\sigma (1) \sigma (2) \tau \circ \sigma (1) \tau
\circ \sigma (2)} \cdot \cdot \cdot R (g)_{\sigma (n-1) \sigma (n) \tau \circ
\sigma (n-1) \tau \circ \sigma (n)} dv_g\,, \leqno (1.2)$$
where $c_n$ is a dimension constant and ${\cal C}_n$ is the set of all
permutations on $\{1, 2,..., n\}$ and $\varepsilon (\tau)$ is the sign of
$\tau \in  {\cal C}_n$. A decomposition of the curvature tensor gives
$$R (g)_{ijkl} = W (g)_{ijkl} + Z (g)_{ijkl} + U (g)_{ijkl}\,, \leqno (1.3)$$
where $W (g)$ is the Weyl curvature tensor and
$$U (g)_{ijkl} = {S (g)\over {n (n - 1)}} (g_{ik} g_{jl} - g_{il}
g_{jk})\,, \leqno (1.4)$$
$$Z (g)_{ijkl} = {1\over {n - 2}} ({\bf z} (g)_{ik} g_{jl} + {\bf z}
(g)_{jl} g_{ik}
- {\bf z} (g)_{il} g_{jk}
- {\bf z} (g)_{jk} g_{il})\,, \leqno (1.5)$$
where ${\bf z } (g)$ is the trace-free Ricci tensor given by
$${\bf z} (g)_{ij} = {\mbox {Ric}}\, (g)_{ij} - {S (g)\over n} g_{ij}\,.
\leqno (1.6)$$
Let $x \in M$ and $\{e_1\,,..., e_n\}$ be an orthonormal basis for the
tangent
space of $M$ above $x$. We have
$$U (g)_{ijkl} = {S (g)\over {n (n - 1)}} (\delta_{ik} \delta_{jl} -
\delta_{il} \delta_{jk}) \ \ \ {\mbox {at}} \ \ x\,.$$
If we apply (1.3), then at the point $x$ we have
\begin{eqnarray*}
(1.7) \ \ \ \ \ \ \ \ \ \ & \ & \sum_{\sigma \in {\cal C}_n} \sum_{\tau
\in {\cal C}_n}
\varepsilon (\sigma) \varepsilon (\tau)
 R (g)_{\sigma (1) \sigma (2) \tau \circ \sigma (1) \tau
\circ \sigma (2)} \cdot \cdot \cdot R (g)_{\sigma (n-1) \sigma (n) \tau \circ
\sigma (n-1) \tau \circ \sigma (n)}\\
& \ \ \ \ \ & \ \ \ \ =
C_o S (g)^{n \over 2} + P (W (g)_{ijkl}, Z
(g)_{ijkl}, U (g)_{ijkl}, g_{ij})\,,
\end{eqnarray*}
where $P$ is a certain polynomial function and $C_o$ is a constant that
depends
on $n$ only. Putting (1.7) into the Gauss-Bonnet formula, we have
\begin{eqnarray*} \chi (M) & = & \int_M C_o S (g)^{n \over 2} dv_g \, +
\, \int_M P (W
(g)_{ijkl}, Z (g)_{ijkl}, U (g)_{ijkl}, g_{ij}) dv_g\\
& = & \int_M C_o S (g')^{n \over 2} dv_{g'} \, + \, \int_M P (W
(g')_{ijkl}, Z
(g')_{ijkl}, U (g')_{ijkl}, g'_{ij})dv_{g'}\,,
\end{eqnarray*}
where $g'$ is another Riemannian metric on $M$. In general, the above
formula is
too complicated to given an effective bounds on $L^{n\over 2}$-norms of
scalar curvatures. We have the following.\\[0.02in]
{\bf Theorem 1.} \ \ {\it Let $(M, h)$ be a compact hyperbolic
n-manifold with $n$ being even.\\[0.02in]
1) \ \ $n = 4$. For any conformally flat metric $g$ on $M$, we have}
$$\int_M |S (g)|^2 dv_g \ge \int_M |S (h)|^2 dv_h\,,$$
{\it and equality holds if and only if $g$ is, up to a positive constant,
isometric to $h$.\\[0.01in]
2) \ \ $n \ge 4$. For any conformally flat metric $g$ on $M$, we have}
$$\int_M |{\mbox {Ric}}\, (g)|^{n\over 2} dv_g \ge c_n \int_M |{\mbox
{Ric}}\, (h)|^{n \over 2} dv_h\,,$$
{\it where $c_n$ is a positive constant that depends on $n$ only.
}\\[0.03in]
{\bf Theorem 2.} \ \ {\it Let $(M, h)$ be a compact hyperbolic
n-manifold with $n$ being even. There exists a positive constant $c'_n$
which depends on
$n$ only, such that for any metric $g$ on $M$ with non-positive sectional
curvature, we have}
$$\int_M |S (g)|^{n\over 2} dv_g \ge c_n' \int_M |S (h)|^{n\over 2}
dv_h\,.$$
\hspace*{0.3in}Besson, Courtois and Gallot [4] have shown that if $(M,
g)$ is a
compact Einstein manifold with negative sectional curvature, then for any
metric
$g'$ in a neighborhood of $g$, we have
$$\int_M |S (g')|^{n \over 2} dv_{g'} \ge \int_M
|S (g)|^{n \over 2} dv_g\,. \leqno (1.8)$$
In the proof of this
result, they have investigated the following.\\[0.02in]
(I) \ \ (1.8)  holds whenever $g'$ is conformal to $g$, i.e.,
if $g' = u^{4\over {n - 2}} g$
for some smooth function $u > 0$ and if $S (g)$ is a {\it negative}
constant, then we have
$$\int_M |S (g')|^{n \over 2} dv_{g'} \ge \int_M
|S (g)|^{n \over 2} dv_g\,.$$
Then they used the second variation formula to investigate the local
behavior of the $L^{n/2}$-norm of $S (g)$. Partially motivated by their
results, we consider the change of
$$\int_M |S (g)|^{n \over 2} dv_g \ \ \ {\mbox {and}}\, \ \ \
\int_M | {\mbox {Ric}}\, (g)
|^{n\over 2} dv_g$$
under Ricci flow and conformal change of metrics when $S (g)$ is a {\it
positive}
constant . The Ricci flow have been considered by Hamilton [11] and by
many authors. It has been proven to be very useful in deforming metrics
into standard metrics, especially
when the original metric is close to a standard metric. For example,
it has been shown in [15, 20] that the Ricci flow starting near a Einstein
metric of negative sectional curvature always converges to it. We obtain
the following
behaviors of $L^{n\over 2}$-norms on curvatures under the Ricci
flow.\\[0.02in]
{\bf Theorem 3.} \ \ {\it Let $(M, g)$ be a compact Riemannian manifold
with $S
(g) < 0$. Let $g_t$ be the Ricci flow starting at $g$, if $S (g_t) \le
0$, then
$${d\over {dt}} \int_M |S (g_t)|^{n \over 2} dv_{g_t} \le 0\,.$$
If we assume that the sectional curvature $K_g$ of $g$ is suitably pinched
$$-1 - \epsilon \le K_g \le -1 + \epsilon$$
for some $\epsilon > 0$, then}
$${d\over {dt}} \int_M | {\mbox {Ric}}\, (g_t)
|^{n\over 2} dv_{g_t} \le 0\,.$$
{\it Under the above conditions, if the Ricci flow converges to a smooth
metric
$g_o$ on $M$, then}
$$\int_M |S (g)|^{n \over 2} dv_{g} \ge
\int_M |S (g_o)|^{n \over 2} dv_{g_o} \ \ \ and$$
$$\int_M | {\mbox {Ric}}\, (g)
|^{n\over 2} dv_{g} \ge \int_M | {\mbox {Ric}}\, (g_o)
|^{n\over 2} dv_{g_o}\,.$$

\hspace*{0.3in}In particular, we provide an alternative proof to (1.8).
In the last section, we consider conformal change of metrics when the
scalar curvature is positive. An interesting question is to what
extend Besson-Courtois-Gallot's
result, viz,
$$\int_M |S (g')|^{n \over 2} dv_{g'} \ge
\int_M |S (g)|^{n \over 2} dv_g$$
if $g'$ is conformal to $g$ and $g$ has constant negative scalar
curvature, holds for
positive scalar curvature.\\[0.02in]
{\bf Theorem 4.}  \ \ {\it Let $(M, g_o)$ be an $n$-manifold with $ b^2 g
\ge {\mbox {Ric}}\, (g) \ge a^2 g$
for some positive numbers $a$ and $b$.
Then for any metric $g = u^{4\over {n - 2}} g_o$, $u > 0$, we have}
$$\int_M |S (g)|^{n \over 2} dv_{g} \ge c_n
\int_M |S (g_o)|^{n \over 2} dv_{g_o}\,,$$
{\it where $c_n$ is a positive constant that depends on $a, b$ and $n$
only.
In general, $c_n < 1$. For the special cases that \\[0.02in]
i) \ \ $g$ is an Einstein metric with positive scalar curvature and
$g = u^{4\over {n - 2}} g_o$, $u > 0$; or\\[0.015in]
ii) \ \ $(M, g)$ is a compact conformally flat manifold with positive
Ricci curvature and
$g_o$ has constant positive sectional curvature;\\[0.02in]
then we have}
$$\int_M |S (g)|^{n \over 2} dv_{g} \ge
\int_M |S (g_o)|^{n \over 2} dv_{g_o}\,.$$

\vspace{0.45in}

{\bf \LARGE 2. \ \ Gauss-Bonnet Formula}

\vspace{0.3in}

Given a compact $n$-manifold $M$ with $n \ge 4$ and a Riemannian metric
$g$ on
$M$, the Weyl conformal curvature tensor can be defined by
$$W (g)_{ijkl} = R (g)_{ijkl} - Z (g)_{ijkl} - U (g)_{ijkl}\,,$$
where $Z (g)$ and $U (g)$ are definied in (1.4) and (1.5), respectively.
Using the fact that
$g^{ij} {\bf z} (g)_{ij} = 0$ and $g^{ik} g^{jl} R (g)_{ijkl} = S (g)$,
it is easy to
show that $g^{ik} g^{jl} W (g)_{ijkl} = 0$ and $g^{ik} W_(g)_{ijkl} = 0$. And
we have
$$|R (g)|^2 = |W (g)|^2 + |Z (g)|^2 + |U (g)|^2\,. \leqno (2.1)$$
A direct calculation shows that
$$| U (g)|^2 = {{2 S (g)^2}\over {n (n - 1)}}\,, \ \ | Z (g)|^2
= {4 \over {(n - 2)}} | {\bf z} (g) |^2 \ \ {\mbox {and}}\, \ \
| {\mbox {Ric}}\, (g)|^2 = | {\bf z} (g)|^2 + {{S (g)^2}\over n}\,.
\leqno (2.2)$$
In dimension four, the Gauss-Bonnet formula takes the form [3]
$$\chi (M) = {1\over {8 \pi^2}} \int_M (|U (g)|^2 - |Z (g)|^2 + |W
(g)|^2) dv_g\,,
\leqno (2.3)$$
where $\chi (M)$ is the Euler characteristic of $M$. Let $h$ be a hyperbolic
metric on $M$, then
$$\chi (M) = {1\over {48 \pi^2}} \int_M S (h)^2 dv_h\,, \leqno (2.4)$$
where $S (h) = - 4 \cdot 3 = - 12$. In dimension bigger than or equal to
four,
a Riemannian metric $g$ is conformally flat if and only if $W (g) \equiv
0$. Then (2.2),
(2.3) and (2.4) show that, if $g$ is any conformally flat metric on $M$,
we have
$$\int_M S (g)^2 dv_g \ge \int_M S (h)^2 dv_h\,.$$
Furthermore, equality holds if and only if ${\bf z} (g) \equiv 0$ and $W
(g) \equiv 0$,
i.e., $(M, g)$ is a hyperbolic metric. By the Mostow rigidity theorem,
$(M, g)$ is
isometric to $(M, h)$ up to a positive constant.\\[0.03in]
{\bf Theorem 2.5.} \ \ {\it Let $(M, h)$ be a compact hyperbolic
$n$-manifold
and $n \ge 4$ being even. For any conformally flat metric $g$ on $M$, we
have}
$$\int_M |{\mbox {Ric}}\, (g)|^{n\over 2} dv_g \ge c_n \int_M |{\mbox
{Ric}}\, (h)|^{n \over 2} dv_g\,,$$
{\it where $c_n$ is a positive constant that depends on $n$
only.}\\[0.03in]
{\it Proof.} \ \ As $n \ge 4$, the metric $g$ is conformally flat if and
only if
$W (g) \equiv 0$. Therefore $R (g) = Z (g) + U (g)$. Apply the
Gauss-Bonnet theorem
we have
$$\chi (M) = C (n) \int_M P (Z (g), U (g)) dv_g\,,$$
where $C (n)$ is a constant that depends on $n$ only and $P$ is a
homogeneous polynomial
of degree $n/2$ in the components of $Z (g)$ and $U (g)$. There exists
positive constants
$C_1, C_1, C_2,..., C_{n\over 2}$, which depend on $n$ only, such that
$$|\chi (M)| \le \int_M ( C_0 |Z (g)|^{n\over 2} + C_1 |Z (g)|^{{n\over
2} - 1} |U (g)| +
C_2 |Z (g)|^{{n\over 2} - 2} |U (g)|^2 + \cdot \cdot \cdot + C_{n\over 2}
|U (g)|^{n\over 2}) dv_g\,.$$
Using (2.2) we have $|{\mbox {Ric}}\, (g)| \ge (n - 2)/\sqrt {4 (n - 2)}
|Z (g)|$ and
$|{\mbox {Ric}}\, (g)| \ge \sqrt {(n - 1)/2} |U (g)|$, we have
$$|\chi (M)| \le C \int_M |{\mbox {Ric}}\, (g)|^{n\over 2} dv_g\,,$$
where $C$ is a constant that depends on $n$ only. For the hyperbolic
metric $h$, we have
$W (h) \equiv 0$, $Z (h) \equiv 0$ and $|{\mbox {Ric}}\, (h)|^2 = S
(h)^2/ n$. The Gauss-Bonnet
theorem gives
$$|\chi (M)| = C' (n) \int_M |S (h)|^{n\over 2} dv_h = C'' (n)
\int_M |{\mbox {Ric}}\, (h)|^{n\over 2} dv_h\,,$$
where $C' (n)$ and $C'' (n)$ are positive constants that depends on $n$
only.
Combining the two formulas we have the result. \ \ \ \ \ {\bf
Q.E.D.}\\[0.03in]
{\bf Theorem 2.6.} \ \ {\it Let $(M, h)$ be a compact hyperbolic
$n$-manifold of
even dimension. There exists a positive constant $c_n$ which
depends on $n$ only, such that for any Riemannian metric $g$ on
$M$ with nonpositive sectional curvature, we have}
$$\int_M S (g)^{n\over 2} dv_g \ge c_n \int_M S (h)^{n \over 2} dv_h\,.$$
{\it Proof.} \ \ By (1.1), we may assume that $n \ge 4$. As $h$ is a
hyperbolic metric, $W (h) \equiv 0$ and $Z (h)
\equiv 0$. The Gauss-Bonnet formula (1.2) gives
$$\chi (M) = \int_M C_o S (h)^{n\over 2} dv_h\,,$$
where $C_o$ is a non-zero constant that depends on $n$ only (its value
can be
found by applying the Gauss-Bonnet formula on $S^n$ and the fact that $\chi
(S^n) = 2$ if $n$ is even). For the Riemannian metric $g$, making use of the
fact that $R (g) = W (g) + Z (g) + U (g)$, the Gauss-Bonnet formula gives
$$\chi (M) = \int_M (C_o S (g)^{n \over 2} + P (W (g)_{ijkl}, Z
(g)_{ijkl}, U (g)_{ijkl}, g_{ij})) dv_g\,,$$
where $P$ is a certain polynomial such that each term contain exactly $n/2$
terms of $W (g)_{jklk}, Z (g)_{ijkl}$ or $U (g)_{ijkl}$. Therefore we have
$$| \chi (M) | \le \int_M C_o | S (g) |^{n \over 2} dv_g + \int_M | P (W
(g)_{ijkl}, Z (g)_{ijkl}, U (g)_{ijkl}, g_{ij})) | dv_g\,.$$
{}From (2.1), $| R(g) | \ge | W (g) |, | R(g) | \ge | Z (g) |$ and
$| R(g) | \ge | U (g) |$, there exists a
positive constant $C_n$ that depends on $n$ only, such that
$$| P (W (g)_{ijkl}, Z (g)_{ijkl}, U (g)_{ijkl}, g_{ij})) | \le C_n | R (g)
|^{n \over 2}\,.$$
Given a point $x \in M$, we choose an orthonormal basis $\{e_1\,,...., e_n\}$
for the tangent space of $M$ above $x$. Let $\sigma_{ij}$ be the
sectional curvature of
the plane spanned by $e_i$ and $e_j$, $i \not= j$, with respect to the
Riemannian metric $g$ on $M$. Assume that $\sigma_{ij} \le 0$. We may also
assume that $\sigma_{12}$ is the minimum of the sectional curvatures at the
point $x$. We have
$$|S (g) |= | \sum_{i,j, i \not= j} \sigma_{ij}| \ge |\sigma_{12}|\,.$$
Let $\sigma ({\bf u}, {\bf v})$ be the sectional curvature of the plane
spanned
by ${\bf u}$ and ${\bf v}$ in the tangent space of $M$ above $x$. Then we
have [7]
\begin{eqnarray*}
R (g)_{ijkl} & = & {1\over 6} \{ 4 [\sigma (e_i + e_l, e_j + e_k) -
\sigma (e_j
+ e_l, e_i + e_k)]\\
& \ & \ - 2 [\sigma (e_i, e_j + e_k) + \sigma (e_j, e_i + e_l) + \sigma (e_k,
e_i + e_l) + \sigma (e_l, e_j + e_k)]\\
& \ & \ + 2 [\sigma (e_i, e_j + e_l) + \sigma (e_j, e_k + e_l) + \sigma (e_k,
e_j + e_l) + \sigma (e_l, e_i + e_k)]\\
& \ & \ + \sigma_{ik} + \sigma_{jl} - \sigma_{il} - \sigma_{jk}\}\,.
\end{eqnarray*}
There exists a positive constant $C'$ which depends on $n$ only, and with
$g_{ij} = \delta_{ij}$, such that we obtain
$$| R
(g) |^2 = \sum_{ijkl} R_{ijkl} R_{ijkl} \le C' (\sigma_{12})^2 \le C' |S (g)
|^2\,,$$ and so
$$| P (W (g)_{ijkl}, Z (g)_{ijkl}, U (g)_{ijkl}, g_{ij})) | \le C_n C'| S (g)
|^{n \over 2}\,.$$
Thus
$$|\chi (M) \le ( |C_o| + C_n C') \int_M | S (g) |^{n \over 2} dv_g\,,$$
or
$$\int_M | S (h) |^{n \over 2} dv_h \le C \int_M | S (g) |^{n \over 2}
dv_g\,,$$
where $C = 1 + C_n C'/|C_o|$ is a positive constant that depends on $n$
only. \ \
\ \ \ {\bf Q.E.D.}\\[0.02in]
{\it Remark:} \ \ From the proof of the above theorem, one can replace the
condition of non-positive sectional curvature by a pinching condition
that the
absolute value of sectional curvature of any 2-plane above a point $x \in M$
is lesser than or equal to $c_n |S (g) (x)|$, a positive constant times the
absolute value of the scalar curvature at that point. Then we have
$$\int_M | S (g) |^{n \over 2} dv_g \ge c' \int_M | S (h) |^{n \over 2}
dv_h\,,$$
where $c'$ is now a constant that depends both on $n$ and $c_n$.\\[0.02in]
{\it Remark:} \ \ It is easy to see that the same result in theorem 2.6
holds for
{\it conformally flat} metrics of nonpositive Ricci curvature.\\[0.02in]
\hspace*{0.3in}The Gauss-Bonnet formula yields the following estimate on
the $L^{n/2}$-norm of scalar curvature.\\[0.03in]
{\bf Lemma 2.7.} \ \ {\it For an even integer $n$ bigger than two,
let $(M, g)$ be a compact $n$-manifold with $\chi (M) \not= 0$. Then
there exist positive
constants $\delta_n$ and $\epsilon_n$ depending on $n$,
which can be a priori choosen, such that, if}
$$\int_M | Z (g) |^{n\over 2} dv_g \le \delta_n \,, \ \ \ \ {\mbox
{and}}\,
\ \ \ \ \int_M | W (g) |^{n\over 2} dv_g \le \epsilon_n$$
{\it then}
$$\int_M |S (g)|^{n\over 2} dv_g \ge c_n\,,$$
{\it where $c_n$ is a positive constant that
depends on $n$ only.}\\[0.03in]
{\it Proof.} \ \ As $R (g) = W (g) + Z (g) + U (g)$. Apply the
Gauss-Bonnet theorem
we have
$$\chi (M) = C (n) \int_M P (W (g), Z (g), U (g)) dv_g\,,$$
where $C (n)$ is a non-zero constant that depends on $n$ only and $P$ is
a homogeneous polynomial
of degree $n/2$ in the components of $W (g)$, $Z (g)$ and $U (g)$. There
exist positive constants
$C'_o, C_o, C_1, C_2,..., C_{n\over 2}$ and $C (n_1, n_2, n_3)$ which
depend on
$n$, $n_1$, $n_2$ and $n_3$
only, such that
\begin{eqnarray*}
(2.8) & \ & |\chi (M)|\\
& \le &\int_M ( C'_0 |U (g)|^{n\over 2} + \sum_{n_1, n_2, n_3} C (n_1,
n_2, n_3)
\int_M |U (g)|^{n_1} |Z (g)|^{n_2} |W (g)|^{n_3} ) dv_g\\
& + & \int_M ( C_o |Z (g)|^{n\over 2} + C_1 |Z (g)|^{{n\over 2} - 1}
|W (g)| +
C_2 |z (g)|^{{n\over 2} - 2} |W (g)|^2 + \cdot \cdot \cdot +
C_{n\over 2} |W (g)|^{n\over 2}) dv_g\,,
\end{eqnarray*}
where $n_1, n_2,$ and $n_3$ are positive integers
such that $n_1 + n_2 + n_3 = n/2$ and $n_1 < n/2$.
For positive numbers $s$, $t$, $p$ and $q$ such that
$$s + t = {n\over 2} \ \ \ {\mbox {and}} \ \ \ {1\over p} + {1\over q} =
1\,,$$
a calculation shows that if $t q = n/2$, then we have $s p = n/2$ as
well. Apply the
H\"older's inequality  to (2.7) (twice to the terms with $n_1, n_2,$ and
$n_3$), we have
\begin{eqnarray*}
(2.9) \ \ & \ & |\chi (M)|\\
 & \le & C'_o \int_M |U (g)|^{n\over 2} dv_g\\
 & + & \sum_{n_1, n_2, n_3} C (n_1, n_2, n_3)
( \int_M |U (g)|^{n\over 2} dv_g)^{1\over {p_{n_1}}}
( \int_M |Z (g)|^{n \over 2} dv_g )^{1\over {r_{n_2}}}
( \int_M |W (g)|^{n\over 2} dv_g )^{1\over {q_{n_3}}}\\
& + & C_o \int_M |Z (g)|^{n\over 2}
dv_g + C_1 ( \int_M |Z (g)|^{n\over 2} dv_g )^{1\over {p_1}}
( \int_M |W (g)|^{n\over 2} dv_g )^{1\over {q_1}} + \cdot \cdot \cdot\\
& + & C_{{n\over 2} - 1}
( \int_M |Z (g)|^{n\over 2} dv_g )^{1\over {p_{{n\over 2} -1}} }
( \int_M |W (g)|^{n\over 2} dv_g )^{1\over {q_{ {n\over 2} - 1}}}
 + C_{n\over 2} \int_M |W (g)|^{n\over 2} dv_g\,,
\end{eqnarray*}
where
$$p_{n_1}, r_{n_2}, q_{n_3}, \ p_1,..., p_{{n\over 2} - 1} \ \ {\mbox
{and}} \ \
q_1,..., q_{{n\over 2} - 1}$$
are positive constants specifed in the H\"older's inequality.
If we choose $\delta_n$ and $\epsilon_n$ small (which depend on $C_o,
C_1,..., C_{n/2}$,
i.e., depend on $n$ only)
such that
$$\int_M |Z (g)|^{n\over 2} dv_g \le \delta_n \ \ \ \ \
{\mbox {and}} \ \ \ \ \ \int_M |W (g)|^{n\over 2} dv_g \le \epsilon_n$$
so that
\begin{eqnarray*}
& \ & \ \ C_o \int_M |Z (g)|^{n\over 2} dv_g
+ C_1 ( \int_M |Z (g)|^{n\over 2} dv_g )^{1\over {p_1}}
( \int_M |W (g)|^{n\over 2} dv_g )^{1\over {q_1}}
+ \cdot \cdot \cdot\\
& \ & \ \ \ + C_{{n\over 2} - 1}
( \int_M |Z (g)|^{n\over 2} dv_g )^{1\over {p_{ {n\over 2} -1}} }
( \int_M |W (g)|^{n\over 2} dv_g )^{1\over {q_{ {n\over 2} - 1}}}
+ C_{n\over 2} \int_M |W (g)|^{n\over 2} dv_g \le {1\over 2}\,,
\end{eqnarray*}
and the fact that
$$| U (g)|^2 = {{2 S (g)^2}\over {n (n - 1)}}\,,$$
(2.9) gives
$$\int_M |S (g)|^{n\over 2} dv_g \ge 2 c_n ( |\chi (M)| - {1\over 2} )\ge
c_n\,,$$
as $\chi (M) \not= 0$ and hence $|\chi (M)| \ge 1$. Here $c_n$ is a
positive constant that depends
on $n$ only. \ \ \ \ \ {\bf Q.E.D.}\\[0.03in]
{\bf Corollary 2.10.} \ \ {\it For an even integer $n$ bigger than two,
let $(M, g)$ be a compact Einstein $n$-manifold with Ric $(g) = \pm (n -
1) g$. Suppose that
$\chi (M) \not= 0$ and}
$$  \int_M | W (g) |^{n\over 2} dv_g \le {1\over {2 C_{n\over 2}}}$$
{\it then Vol $(M, g) \ge {c'}_n$, where $C_{n\over 2}$
is the same constant as in (2.9) and
${c'}_n$ is a positive constant that
depends on $n$ only.}\\[0.03in]
{\it Proof.} \ \ As $(M, g)$ is an Einstein manifold, we have $Z (g) =
0$. Therefore
in (2.9), the terms involving $Z (g)$ vanish and we just need
$$ \int_M | W (g) |^{n\over 2} dv_g \le {1\over {2 C_{n\over 2}}}$$
to conclude that
$$\int_M |S (g)|^{n\over 2} dv_g \ge c_n\,.$$
Using the fact that
$$|S (g)| = n (n - 1)$$
for an Einstein manifold with Ric $(g) = \pm (n - 1) g$, we obtain the
result.
 \ \ \ \ \ {\bf Q.E.D.}\\[0.03in]
{\bf Corollary 2.11.} \ \ {\it For an even integer $n$ bigger than two,
let $(M, g)$ be a compact Einstein $n$-manifold with Ric $(g) = (n - 1) g$
and $\chi (M) \not= 0$. Then there exists a positive number
$\varepsilon_n$,
which depends on $n$ only, such that if}
$$\int_M | W (g) |^{n\over 2} dv_g \le \varepsilon_n\,,$$
{\it then $g$ has constant positive sectional curvature. In case $n = 4$,
we can
drop the assumption that $\chi (M) \not= 0$.}\\[0.03in]
{\it Proof.} \ \ If we take $\varepsilon_n < 1/(2 C_{n\over 2})$, then
corollary
(2.9) shows that
Vol $(M, g) \ge c'_n$ for some positive constant $c'_n$ that depends on
$n$ only.
A result in [17] shows that there exists a positive constant $c"_n$ which
depends on
$n$ only, such that if
$$\int_M | W (g) |^{n\over 2} dv_g \le c"_n {\mbox {Vol}}\, (M, g)\,,$$
then $g$ is a metric of constant sectional curvature one. We can take
$\epsilon_n =  \min \{ c'_n c"_n, 1/(2 C_{n\over 2}) \}$.
If $n = 4$, then the Gauss-Bonnet formula for an Einstein metric has the form
$$\chi (M) = {1\over {8 \pi^2}} \int_M (|U (g)|^2 + |W (g)|^2) dv_g\,.$$
It follows that $\chi (M) \not= 0$ if Ric $(g) = (n - 1) g$.
\ \ \ \ \ {\bf Q.E.D.}\\[0.3in]
{\it Remark:} \ \ Similar pinching results are obtained in [17] and
[9].\\[0.03in]
\hspace*{0.3in}So far the best estimates we get are on K\"ahler-Einstein
manifolds
of negative scalar curvature ($n = 4$).

{\bf Theorem 2.12.} \ \ {\it Let $\eta$ be a
K\"ahler-Einstein metric of negative scalar curvature on a compact complex
surface $M$. Then for any Riemannian metric $g$ on $M$, we have}
$$\int_M |S (g)|^2 dv_{g} \ge \int_M |S (\eta)|^2 dv_\eta\,.$$

{\it Proof.} \ \  This follows rather directly from Seiberg-Witten theory
[19, 13]. The main
point is this, for the K\"ahler-Einstein surface $(M, \eta)$, the first Chern
class is given by $[S (\eta) \omega / 8 \pi]$, where $\omega$ is the K\"ahler
form of
$(M, \eta)$. Hence
$$c_1 (L)^2 = 2 \chi (M) + 3 \tau (M) = {1\over {32 \pi^2}} \int_M |S
(\eta)|^2
dv_\eta\,.$$ Here $\chi (M)$ and $\tau (M)$ are the Euler characteristic and
signature of
$M$, respectively. It follows from Seiberg-Witten theory that for any
Riemannian
metric $g$ on $M$, there exists a solution to the Seiberg-Witten
equations [19]
and we have [13]
$$c_1 (L)^2 \le {1\over {32 \pi^2}} \int_M |S (g)|^2 dv_g\,.$$ That is,
$$\int_M |S (g)|^{2} dv_{g} \ge \int_M |S (\eta)|^{2} dv_\eta\,.
$$
This completes the proof. \ \ \ \ \ {\bf Q.E.D.}

\vspace{0.45in}

{\bf \LARGE \S 3. \ \ Ricci Curvature Flow}

\vspace{0.3in}

Let $(M, g_o)$ be a compact Riemannian manifold. In this section we
consider the
Ricci curvature flow:
$${{\partial g}\over {\partial t}} = - 2 {\bf z} (g) - {{2 \delta S (g)
}\over n} g\,, \ \ \ \ \ g (0) = g_o\,,  \leqno (3.1)$$
where ${\bf z} (g) = {\mbox {Ric}}\, (g) - [S (g)/n] g$ is the trace free
Ricci
tensor as in section 1 and
$$\delta S (g) = S (g) - {{\int_M S (g) dv_g}\over {\int_M dv_g}}\,.$$
The Ricci curvature flow has been studied extensively by Hamilton,
Huisken, Margerin, Nishikawa, Shi, Ye, and many others in respect to the
questions of long time existence and convergence, we refer to [17] for a
comprehensive
references. It has been shown that if $(M, h)$ is a compact
Einstein manifold of strictly negative sectional curvature, then there
exists an open neighborhood of $h$ in the space of smooth metrics with
$C^\infty$-norm
such that each metric $g_o$ in that open neighborhood converges to $h$
under the
Ricci curvature flow (3.1)  [15, 20]. Furthermore, we can choose an open
neighborhood such that the Ricci curvature remain negative during the
Ricci curvature flow.\\[0.03in]
{\bf Lemma 3.2.} \ \ {\it For $n \ge 4$, let $M$ be a compact
$n$-manifold. Let
$g$ be a solution to the Ricci curvature flow equation (3.1) on the time
interval
$(0, t')$, where $t'$ may equal to infinity. Assume that $\lim_{t \to t'}
g =
g'$ is a smooth Riemannian metric on $M$. If $S (g) < 0$ for $t \in (0,
t')$, then}
$${d\over {dt}} \int_M |S (g)|^{n\over 2} dv_g \le 0 \ \ \ \ {\mbox {for \
all}} \ \ \ t \in (0, t')\,.$$
{\it Hence}
$$\int_M |S (g_o)|^{n\over 2} dv_{g_o} \ge \int_M |S (g')|^{n\over 2}
dv_{g'}\,.$$
{\it Proof.} \ \ From (3.1) we have [11,16]
$${{d S (g)}\over {dt}} = \Delta S (g) + 2 |{\bf z}|^2 + {{2 \delta S (g)
}\over n} S (g)\,. \leqno (3.3)$$
As ${\bf z} (g)$ is trace-free, we have
$$(dv_g)' = {1\over 2} {\mbox {tr}}_g ({{dg}\over {dt}}) dv_g  = - \delta S
(g)\,. \leqno (3.4)$$
Therefore
\begin{eqnarray*}
& \ & {d \over {dt}} \int_M |S (g)|^{n\over 2} dv_g\\
 & = & \int_M ({n\over 2}
|S (g)|^{{n\over 2} - 1} {d \over {dt}} |S (g)| dv_g + \int_M |S
(g)|^{n\over 2}
 (dv_g)'\\
& = & \int_M ({n\over 2}
|S (g)|^{{n\over 2} - 1} (- {d \over {dt}} S (g)) dv_g - \int_M |S
(g)|^{n\over
2} (\delta S (g)) dv_g \ \ \ \ ({\mbox {as}} \ \ S (g) < 0)\\
& = & \int_M ({n\over 2}
|S (g)|^{{n\over 2} - 1} (- \Delta S (g) - 2 |{\bf z}|^2 - {{2 \delta S (g)
}\over n} S (g)) dv_g - \int_M |S (g)|^{n\over
2} (\delta S (g)) dv_g\\
& = & - \int_M ({n\over 2}) ({n\over 2} - 1) |S (g)|^{{n\over 2} - 2}
|\bigtriangledown |S (g)||^2 - 2 \int_M {n\over 2}
|S (g)|^{{n\over 2} - 1} |{\bf z}|^2 dv_g \le 0\,,
\end{eqnarray*}
as $-\Delta S (g) = \Delta | S (g)|\,.$ \ \ \ \ \ {\bf Q.E.D.}\\[0.05in]
{\bf Theorem 3.5.} [4] \ \ {\it For $n \ge 4$, let $(M, h)$ be a compact
Einstein $n$-manifold of strictly negative sectional curvature. Then
there exists an open neighborhood of $h$ in the space of smooth metrics
on $M$
with $C^\infty$-norm such that for any metric $g$ in the open
neighborhood,}
$$\int_M |S (g)|^{n\over 2} dv_{g} \ge \int_M |S (h)|^{n\over 2}
dv_{h}\,.$$
{\it Proof.} \ \ The existence of such an open neighborhood of $h$ for which
the the Ricci curvature flow (3.1) converges to $h$ is shown in [17].
Furthermore, we may choose the open neighborhood such that
the scalar curvature remains negative
during the Ricci curvature flow. Then we can apply lemma 3.2.\\[0.05in]
{\bf Theorem 3.6.} \ \ {\it For $n \ge 4$, let $(M, h)$ be a compact
hyperbolic $n$-manifold. Then
there exists an open neighborhood of $h$ in the space of smooth metrics on
$M$ with $C^\infty$-norm such that for any metric $g_o$ in the open
neighborhood, if $g$ is a
solution to the Ricci curvature flow (3.1) with initial condition $g_o$,
then}
$${d\over {dt}} \int_M |{\mbox {Ric}}\, (g)|^{n\over 2} dv_g \le 0\,.$$
{\it Proof.} \ \ As $| {\mbox {Ric}}\, (g)|^2 = |{\bf z} (g)|^2 + S
(g)^2/ n$,
we have
\begin{eqnarray*}
{d\over {dt}} (| {\mbox {Ric}}\, (g)|^{n \over 2}) & = & {d\over {dt}} (|
{\mbox {Ric}}\,
(g)|^2)^{n \over 4})\\
& = & {n\over 4} (| {\mbox {Ric}}\, (g)|^2)^{({n \over 4} - 1)}
{d\over {dt}} | {\mbox {Ric}}\, (g)|^2\\
& = & {n\over 4} (| {\mbox {Ric}}\, (g)|^2)^{({n \over 4} - 1)}
{d\over {dt}} (|{\bf z} (g)|^2 + {{S (g)^2}\over
n})\,.
\end{eqnarray*}
We have [17]
$${{\partial}\over {\partial t}} |{\bf z} (g)|^2 = \Delta |{\bf z} (g)|^2
-
2 | \bigtriangledown {\bf z} (g) |^2 + 4 Rm ({\bf z} (g)) \cdot {\bf z}
(g) +
{4\over n} \delta S (g) |{\bf z} (g)|^2\,,$$
where $Rm ({\bf z} (g)) \cdot {\bf z} (g) = g^{ii'} g^{jj'} g^{kk'}
g^{ll'}
R (g)_{ijkl} {\bf z} (g)_{i'k'} {\bf z} (g)_{j'l'} $. From (3.3) we have
\begin{eqnarray*}
{{\partial}\over {\partial t}} |S|^2 & = & 2 S (g) \Delta S (g) + 4 S (g)
|{\bf z} (g)|^2 + {4\over n} \delta S (g) S (g)^2\\
& = & \Delta |S (g)|^2 - 2 | \bigtriangledown |S (g)||^2
+ 4 S (g)
|{\bf z} (g)|^2 + {4\over n} \delta S (g) S (g)^2\,,
\end{eqnarray*}
as $S (g) < 0$ and $\Delta u^2 = 2 u \Delta u + 2 |\bigtriangledown u|^2$.
Therefore
\begin{eqnarray*}
& \ & \ {d\over {dt}} \int_M |{\mbox {Ric}}\, (g)|^{n\over 2} dv_g\\
 & = & \int_M {n\over 4} [| {\mbox {Ric}}\, (g)|^2)^{({n \over 4} - 1)}
(\Delta |{\bf z} (g)|^2 - 2 | \bigtriangledown {\bf z} (g) |^2 + 4 Rm
({\bf z}
(g)) \cdot {\bf z} (g) +  {4\over n} \delta S (g) |{\bf z} (g)|^2\\
& \ & \ \ \ \  + {1\over n} (\Delta |S (g)|^2 - 2 | \bigtriangledown |S
(g)||^2)
 + {4\over n} S (g)  |{\bf z} (g)|^2 + {4\over n} \delta S (g)
{{S (g)^2}\over n} ] dv_g\\
& \ & \ \ \ \  - \int_M  |{\mbox {Ric}}\, (g)|^{n\over 2} \delta S (g) dv_g\\
& = & \int_M {n\over 4} ({n \over 4} - 1) (| {\mbox {Ric}}\, (g)|^2)^{({n
\over 4} - 2)} (
-  | \bigtriangledown {\mbox {Ric}}\, (g) |^2 + 4 Rm ({\bf z}
(g)) \cdot {\bf z} (g)\\
& \ & \ \ \ \ - {2\over n} |\bigtriangledown |S (g)| |^2 + {4\over n} S
(g)  |{\bf
z} (g)|^2) dv_g\,,
\end{eqnarray*}
as $\Delta |{\bf z} (g)|^2 + {1\over n} \Delta |S (g)|^2 = \Delta |{\mbox
{Ric}}\, (g)|^2$.
Therefore ${d\over {dt}} \int_M |{\mbox {Ric}}\, (g)|^{n\over 2} dv_g \le
0$ if we
can show that
$$Rm ({\bf z} (g)) \cdot {\bf z} (g) + {1\over n} S (g)  |{\bf z}
(g)|^2 \le 0\,.$$

{\bf Lemma 3.7.} \ \ {\it There exists a positive constant $\epsilon$ which
depends on $n$ only ($n \ge 4$) such that if $(M, g)$ is a compact Riemannian
$n$-manifold with sectional curvature $K$ satisfying $-1 - \epsilon \le K
\le -1 +
\epsilon$, then}
$$n Rm ({\bf z} (g)) \cdot {\bf z} (g) + S (g)  |{\bf z}
(g)|^2 \le 0\,.$$
{\it Proof.} \ \ We show the case $n = 4$ first. Let $x \in M$. Choose an
orthonormal basis $\{e_1, e_2, e_3, e_4 \}$ for the tangent space above
$x$ such
that, at the point $x$,   $$g_{ij} = \delta_{ij} \ \ \ {\mbox {and}} \ \
\ {\bf z}
(g)_{ij} = \lambda_i \delta_{ij} \ \ \ {\mbox {for}} \ \ 1 \le i, j \le
4\,.$$
Let $\sigma_{ij}$ be the sectional curvature of the plane spanned by
$e_i$ and
$e_j$. Then, at the point $x \in M$,
\begin{eqnarray*}
Rm ({\bf z} (g)) \cdot {\bf z} (g) & = & R (g)_{ijkl} {\bf z} (g)_{i'k'}
{\bf z} (g)_{j'l'} g^{ii'} g^{jj'} g^{kk'} g^{ll'}\\
& = & \sum_{i \not= j} R (g)_{ijij} {\bf z} (g)_{ii} {\bf z} (g)_{jj}\\
& = & \sum_{i \not= j} \sigma_{ij} \lambda_i \lambda_j\,.
\end{eqnarray*}
Therefore
\begin{eqnarray*}
4 Rm ({\bf z} (g)) \cdot {\bf z} (g) + S (g)  |{\bf z}
(g)|^2) & = & 4 \sum_{i \not= j} \sigma_{ij} \lambda_i \lambda_j
+ \sum_{i \not= j} \sigma_{ij} (\lambda_1^2 +
\lambda_2^2 + \lambda_3^2 + \lambda_4^2)\\
& = & \sum_{i \not= j} \sigma_{ij} (4\lambda_i \lambda_j
+ \lambda_1^2 +
\lambda_2^2 + \lambda_3^2 + \lambda_4^2)\,.
\end{eqnarray*}
We need to show that
$$\sum_{i \not= j} \sigma_{ij} (4\lambda_i \lambda_j
+ \lambda_1^2 +
\lambda_2^2 + \lambda_3^2 + \lambda_4^2) \le 0\,. \leqno (3.8)$$
Assume that $-1 \le \sigma_{ij} \le -1 + \epsilon$ for $1 \le i, j \le 4$.
Then
\begin{eqnarray*}
(3.9) \ \ \ \ \ \ & \ &\sigma_{12} (4\lambda_1 \lambda_2
+ \lambda_1^2 +
\lambda_2^2 + \lambda_3^2 + \lambda_4^2)
+ \sigma_{34} (4\lambda_3 \lambda_4
+ \lambda_1^2 +
\lambda_2^2 + \lambda_3^2 + \lambda_4^
2)\\
& = & -2[(\lambda_1 + \lambda_2)^2 + (\lambda_3 + \lambda_4)^2] + O
(\epsilon)
[4 (\lambda_1 \lambda_2 + \lambda_3 \lambda_4)
+ 2 (\lambda_1^2 +
\lambda_2^2 + \lambda_3^2 + \lambda_4^2)]\,.
\end{eqnarray*}
And
\begin{eqnarray*}
(3.10) \ \ \ \ \ \ \ \ & \ &\sigma_{13} (4\lambda_1 \lambda_3
+ \lambda_1^2 +
\lambda_2^2 + \lambda_3^2 + \lambda_4^2)
+ \sigma_{14} (4\lambda_1 \lambda_4
+ \lambda_1^2 +
\lambda_2^2 + \lambda_3^2 + \lambda_4^
2)\\
& \ & \ \ \ \ + \sigma_{23} (4\lambda_2 \lambda_3
+ \lambda_1^2 +
\lambda_2^2 + \lambda_3^2 + \lambda_4^2)
+ \sigma_{24} (4\lambda_2 \lambda_4
+ \lambda_1^2 +
\lambda_2^2 + \lambda_3^2 + \lambda_4^
2)\\
& = & - [ 2 \lambda_1 \lambda_3 + (\lambda_1 + \lambda_3)^2 + \lambda_2^2 +
\lambda_4^2 + 2\lambda_1 \lambda_4 + (\lambda_1 + \lambda_4)^2 +
\lambda_2^2 +
\lambda_3^2\\
& \ & \ \ \ \ + 2\lambda_2 \lambda_3 + (\lambda_2 + \lambda_3)^2 +
\lambda_1^2 +
\lambda_4^2 + 2\lambda_2 \lambda_4 (\lambda_2 + \lambda_4)^2 +
\lambda_1^2 +
\lambda_3^2]\\
& \ & \ \ \ \ + O (\epsilon) [4( \lambda_1 \lambda_3 + \lambda_1
\lambda_4
+ \lambda_2 \lambda_3 + \lambda_2 \lambda_4
+ \lambda_1^2 +
\lambda_2^2 + \lambda_3^2 + \lambda_4^
2)]\,.
\end{eqnarray*}
Since
$$-[(\lambda_1 + \lambda_2)^2 + (\lambda_3 + \lambda_4)^2 +
2 (\lambda_1 \lambda_3 + \lambda_1 \lambda_4
+ \lambda_2 \lambda_3 + \lambda_2 \lambda_4)] = - [(\lambda_1 +
\lambda_2) +
(\lambda_3 + \lambda_4)]^2\,,$$
we add (3.9) and (3.10) together to obtain
\begin{eqnarray*}
(3.11) \ \ \ \ \ & \ &\sigma_{12} (4\lambda_1 \lambda_2
+ \lambda_1^2 +
\lambda_2^2 + \lambda_3^2 + \lambda_4^2)
+ \sigma_{34} (4\lambda_3 \lambda_4
+ \lambda_1^2 +
\lambda_2^2 + \lambda_3^2 + \lambda_4^
2)\\
& \ & \ \ \ \ + \sigma_{13} (4\lambda_1 \lambda_3
+ \lambda_1^2 +
\lambda_2^2 + \lambda_3^2 + \lambda_4^2)
+ \sigma_{14} (4\lambda_1 \lambda_4
+ \lambda_1^2 +
\lambda_2^2 + \lambda_3^2 + \lambda_4^
2)\\
& \ & \ \ \ \ + \sigma_{23} (4\lambda_2 \lambda_3
+ \lambda_1^2 +
\lambda_2^2 + \lambda_3^2 + \lambda_4^2)
+ \sigma_{24} (4\lambda_2 \lambda_4
+ \lambda_1^2 +
\lambda_2^2 + \lambda_3^2 + \lambda_4^
2)\\
& = & - [(\lambda_1 + \lambda_2 + \lambda_3 + \lambda_4)^2
+ (\lambda_1 + \lambda_2)^2 +
(\lambda_3 + \lambda_4)^2 + (\lambda_1 + \lambda_3)^2\\
& \ & \ \ \ \ + (\lambda_1 +
\lambda_4)^2 + (\lambda_2 + \lambda_4)^2 + (\lambda_2 + \lambda_4)^2 +
2 (\lambda_1^2
+ \lambda_2^2 + \lambda_3^2 + \lambda_4^2)]\\
& \ & \ \ \ \ + O (\epsilon) [4( \lambda_1 \lambda_2 + \lambda_3
\lambda_4 +
\lambda_1 \lambda_3 + \lambda_1 \lambda_4  + \lambda_2 \lambda_3 + \lambda_2
\lambda_4)\\
& \ & \ \ \ \ + 6 (\lambda_1^2 +
\lambda_2^2 + \lambda_3^2 + \lambda_4^
2)] \le 0\,.
\end{eqnarray*}
The last inequality holds if we choose $\epsilon$ to be small, as
the term $(\lambda_1^2 + \lambda_2^2 + \lambda_3^2 + \lambda_4^2)$ will
dominate all the terms with $\epsilon$. We can explicitly choose
$\epsilon = 1/4$. As
$\sigma_{ij} = \sigma_{ji}$, the remaining six terms in (3.8) is in fact the
same as in (3.11). Hence $$4 Rm ({\bf z} (g)) \cdot {\bf z} (g) + S
(g)  |{\bf z} (g)|^2 \le 0\,.$$
For $n > 4$, the proof is similar but more complicated. Choose an
orthonormal basis for the tangent space above $x \in M$
such that ${\bf z} (g)_{ij} = \lambda_i \delta_{ij}$ for $1 \le i, j \le
n$.
We need to show that $$\sum_{i < j\,, 1 \le i, j \le n} \sigma_{ij} (n
\lambda_i
\lambda_j + \lambda_1^2 + \cdot \cdot \cdot + \lambda_{n - 1}^2 +
\lambda_n^2)
\le 0\,.$$ By induction, we may assume that there exists a positive
number $c_{n
- 1}$ such that
\begin{eqnarray*}
& \ & \ \sum_{i < j\,, 1 \le i, j \le n - 1} \sigma_{ij} [(n - 1) \lambda_i
\lambda_j + \lambda_1^2 + \cdot \cdot \cdot + \lambda_{n - 1}^2]\\
& \ & \le - c_{n - 1} (\lambda_1^2 + \cdot \cdot \cdot + \lambda_{n -
1}^2) +
O (\epsilon) (\sum_{i < j\,, 1 \le i, j \le n - 1}\lambda_i
\lambda_j + \lambda_1^2 + \cdot \cdot \cdot + \lambda_{n - 1}^2)
\end{eqnarray*}
Then
\begin{eqnarray*}
& \ & \ \sum_{i < j\,, 1 \le i, j \le n - 1} \sigma_{ij} [n\lambda_i
\lambda_j + \lambda_1^2 + \cdot \cdot \cdot + \lambda_{n - 1}^2 +
\lambda_n^2]\\
& \ & \le - c_{n - 1} (\lambda_1^2 + \cdot \cdot \cdot + \lambda_{n -
1}^2) -
\sum_{i < j\,, 1 \le i, j \le n - 1}\lambda_i \lambda_j - {{(n - 1) (n -
2)}\over 2} \lambda_n^2\\
& \ & \ \ \ \ + O (\epsilon) (\sum_{i < j\,, 1 \le i, j \le n - 1}\lambda_i
\lambda_j + \lambda_1^2 + \cdot \cdot \cdot + \lambda_{n - 1}^2 +
\lambda_n^2)
\end{eqnarray*}
In the sum $\sum_{i < j\,, 1 \le i, j \le n - 1}\lambda_i \lambda_j$, each
$\lambda_i$ appears $n - 2$ times for $1 \le i \le n - 1$. We have
\begin{eqnarray*}
& \ & \sum_{i < j\,, 1 \le i, j \le n} \sigma_{ij} (n \lambda_i \lambda_j +
\lambda_1^2 + \cdot \cdot \cdot + \lambda_{n - 1}^2 + \lambda_n^2)\\
& \le & - c_{n - 1} (\lambda_1^2 + \cdot \cdot \cdot + \lambda_{n - 1}^2)
-
\sum_{i < j\,, 1 \le i, j \le n - 1}\lambda_i \lambda_j - {{(n - 1) (n -
2)}\over 2} \lambda_n^2\\
& \ & \ \ \ \ - (n \lambda_1 \lambda_n
+ \lambda_1^2 + \cdot \cdot \cdot + \lambda_{n - 1}^2 + \lambda_n^2)\\
& \ &  \ \ \ \ \ \ \ \ \ \ \cdot\\
& \ &  \ \ \ \ \ \ \ \ \ \ \cdot\\
& \ &  \ \ \ \ \ \ \ \ \ \ \cdot\\
& \ & \ \ \ \ - (n \lambda_{n - 1} \lambda_n
+ \lambda_1^2 + \cdot \cdot \cdot + \lambda_{n - 1}^2 + \lambda_n^2)\\
& \ & \ \ \ \ + O (\epsilon) (\sum_{i < j\,, 1 \le i, j \le n}\lambda_i
\lambda_j + \lambda_1^2 + \cdot \cdot \cdot + \lambda_{n - 1}^2 +
\lambda_n^2)\\
& = &
- c_{n - 1} (\lambda_1^2 + \cdot \cdot \cdot + \lambda_{n - 1}^2) -
\sum_{i < j\,, 1 \le i, j \le n - 1}\lambda_i \lambda_j - n (\lambda_1
\lambda_n +
\cdot \cdot \cdot + \lambda_{n - 1} \lambda_n)\\
& \ & \ \ \ \ \ - (n - 1) (\lambda_1^2 + \cdot \cdot \cdot +
\lambda_{n - 1}^2)  - [
{{(n - 1) (n -
2)}\over 2} + (n - 1)]\lambda_n^2\\
& \ & \ \ \ \ + O (\epsilon) (\sum_{i < j\,, 1 \le i, j \le n}\lambda_i
\lambda_j + \lambda_1^2 + \cdot \cdot \cdot + \lambda_{n - 1}^2 +
\lambda_n^2)\\
& = & - c_{n - 1} (\lambda_1^2 + \cdot \cdot \cdot + \lambda_{n - 1}^2)
- \sum_{i < j\,, 1 \le i, j \le n - 1}
({1\over 2} \lambda_i + \lambda_i \lambda_j + {1\over 2} \lambda_j)\\
& \ & \ \ \ \ - {n\over 2} (\lambda_1^2 + \cdot \cdot \cdot +
\lambda_{n - 1}^2) - n (\lambda_1 \lambda_n + \cdot \cdot \cdot +
\lambda_{n -1}
\lambda_n) - ({n\over 2})(n - 1) \lambda_n^2\\
& \ & \ \ \ \ + O (\epsilon) (\sum_{i < j\,, 1 \le i, j \le n}\lambda_i
\lambda_j + \lambda_1^2 + \cdot \cdot \cdot + \lambda_{n - 1}^2 +
\lambda_n^2)\\
& = & - c_{n - 1} (\lambda_1^2 + \cdot \cdot \cdot + \lambda_{n - 1}^2)
- {1\over 2} (\lambda_1 + \lambda_2)^2  - (\lambda_3 + \lambda_n)^2\\
& \ & \ \ \ \ - \sum_{i < j\,, 1 \le  i, j \le n - 1\,, (i, j) \not= (1, 2)}
({1\over 2} \lambda_i^2 + \lambda_i \lambda_j + {1\over 2} \lambda_j^2)\\
& \ & \ \ \ \ - ({n\over 2} \lambda_1^2 + {n\over 2} \lambda_2^2 +
({n\over 2} - 1) \lambda_3^2 {n\over 2} + \lambda_4^2 \cdot \cdot \cdot +
{n\over 2} \lambda_{n - 1}^2)\\
& \ & \ \ \ \ - (n \lambda_1 \lambda_n + n \lambda_2 \lambda_n
 + (n - 2) \lambda_3 \lambda_n + n \lambda_4 \lambda_n + \cdot \cdot
\cdot + n \lambda_{n -1}
\lambda_n) - [({n\over 2})(n - 1) - 1] \lambda_n^2\\
& \ & \ \ \ \ + O (\epsilon) (\sum_{i < j\,, 1 \le i, j \le n}\lambda_i
\lambda_j + \lambda_1^2 + \cdot \cdot \cdot + \lambda_{n - 1}^2 +
\lambda_n^2)\,.\\
\end{eqnarray*}
We have
$${1\over 2} (\lambda_1 + \lambda_2)^2  + (\lambda_3 + \lambda_n)^2 +
{\sqrt 2} (\lambda_1 \lambda_3
+ \lambda_1 \lambda_n + \lambda_2 \lambda_3 + \lambda_2 \lambda_n)
= ({1\over {\sqrt 2}} \lambda_1 + {1\over {\sqrt 2}} \lambda_2 +
\lambda_3 + \lambda_n)^2\,.$$
Therefore
\begin{eqnarray*}
& \ & \sum_{i < j\,, 1 \le i, j \le n} \sigma_{ij} (n \lambda_i \lambda_j +
\lambda_1^2 + \cdot \cdot \cdot + \lambda_{n - 1}^2 + \lambda_n^2)\\
& \le & - c_{n - 1} (\lambda_1^2 + \cdot \cdot \cdot + \lambda_{n - 1}^2)
-
({1\over {\sqrt 2}} \lambda_1 + {1\over {\sqrt 2}} \lambda_2 + \lambda_3
+ \lambda_n)^2\\
& \ & \ \ \ \ - \sum_{i < j\,, 1 \le  i, j \le n - 1, (i, j) \not= (1,
2), (1,3), (2, 3)}
{1\over 2} (\lambda_i + \lambda_j)^2 -
({1\over 2} \lambda_1^2 - ({\sqrt 2} - 1) \lambda_1 \lambda_3 + {1\over
2} \lambda_3^2)\\
& \ & \ \ \ \ -({1\over 2} \lambda_2^2 - ({\sqrt 2} - 1) \lambda_2
\lambda_3
+ {1\over 2} \lambda_3^2)\\
& \ & \ \ \ \ - ({n\over 2} \lambda_1^2 + {n\over 2} \lambda_2^2 +
({n\over 2} - 1) \lambda_3^2 + {n\over 2} \lambda_4^2 + \cdot \cdot \cdot +
{n\over 2} \lambda_{n - 1}^2)\\
& \ & \ \ \ \ - [(n - {\sqrt 2}) \lambda_1 \lambda_n
+ (n - {\sqrt 2}) \lambda_2 \lambda_n
 (n - 2) \lambda_3 \lambda_n + n \lambda_4 \lambda_n + \cdot \cdot \cdot
+ n \lambda_{n -1}
\lambda_n)\\
& \ & \ \ \ \ - [({n\over 2})(n - 1) - 1] \lambda_n^2
+ O (\epsilon) (\sum_{i < j\,, 1 \le i, j \le n}\lambda_i
\lambda_j + \lambda_1^2 + \cdot \cdot \cdot + \lambda_{n - 1}^2 +
\lambda_n^2)\\
& \le & - c_n (\lambda_1^2 + \cdot \cdot \cdot + \lambda_{n - 1}^2 +
\lambda_n^2) -
({1\over {\sqrt 2}} \lambda_1 + {1\over {\sqrt 2}} \lambda_2 + \lambda_3
+ \lambda_4)^2\\
& \ & \ \ \ \ - \sum_{i < j\,, 1 \le  i, j \le n - 1, (i, j) \not= (1,
2), (1,3), (2, 3)}
({1\over 2} (\lambda_i + \lambda_j)^2 - {{ {\sqrt 2} - 1}\over 2}
[(\lambda_1 - \lambda_3)^2
 + (\lambda_2 + \lambda_3)^2]\\
& \ & \ \ \ \ {n\over 2} [(\lambda_3 + \lambda_n)^2 + \cdot \cdot \cdot
(\lambda_{n - 1} + \lambda_n)^2]
+ {{n - {\sqrt 2}}\over 2} [(\lambda_1 + \lambda_n)^2 + (\lambda_2 +
\lambda_n)^2]\\
& \ & \ \ \ \ + O (\epsilon) (\sum_{i < j\,, 1 \le i, j \le n}\lambda_i
\lambda_j + \lambda_1^2 + \cdot \cdot \cdot + \lambda_{n - 1}^2 +
\lambda_n^2)
\,,
\end{eqnarray*}
where in the last inequality $c_n$ is a positive constant. Therefore
\begin{eqnarray*}
& \ & \ \sum_{i < j\,, 1 \le i, j \le n} \sigma_{ij} [n \lambda_i
\lambda_j + \lambda_1^2 + \cdot \cdot \cdot + \lambda_{n - 1}^2 +
\lambda_n^2]\\
& \ & \le - c_n (\lambda_1^2 + \cdot \cdot \cdot + \lambda_{n - 1}^2 +
\lambda_n) +
O (\epsilon) (\sum_{i < j\,, 1 \le i, j \le n}\lambda_i
\lambda_j + \lambda_1^2 + \cdot \cdot \cdot + \lambda_{n - 1}^2 +
\lambda_n^2)\,.
\end{eqnarray*}
Hence we can choose $\epsilon$ small such that
$$\sum_{i < j\,, 1 \le i, j \le n} \sigma_{ij} (n \lambda_i \lambda_j +
\lambda_1^2 + \cdot \cdot \cdot + \lambda_{n - 1}^2 + \lambda_n^2) \le 0\,.$$
By induction, we have finished the proof for all $n \ge 4$.
\ \ \ \ \ {\bf Q.E.D.}\\[0.02in]

{\it Proof of Theorem 3.6 Continued.} \ \ We may choose an open
neighborhood of
$h$ such that the sectional curvatures of all the metrics in the open
neighborhood is sufficiently pinched. As shown in [17], curvature
pinching is
preserved during the Ricci curvature flow. Therefore we can apply Lemma
3.7 to
finish the proof.\\[0.02in]
{\it Remark:} \ \ We may apply lemma 3.7 to show theorem 3 in the
introduction.

\vspace{0.3in}

{\bf \LARGE \S 4. \ \ Conformal Changes of Metrics}

\vspace{0.3in}

We begin with the following lemma (cf. [4]), which says that among all
conformal metrics, the
ones with constant nonpositive scalar curvatures have minimal $L^{n\over
2}$-norms of scalar curvatures.
The result has been proved in [4]. For the sake of completeness we
present a proof here,
using a different scalar curvature equation.\\[0.03in]
{\bf Lemma 4.1.} \ \ {\it Let $M$ be a compact $n$-manifold with $n \ge
3$ and $g$ be
a Riemannian metric on $M$ with constant nonpositive scalar curvature.
Then for any metric $g'$ that
is conformal to $g$, we have}
$$\int_M |S (g')|^{n \over 2} dv_{g'} \ge \int_M |S (g)|^{n \over 2}
dv_g\,,$$
{\it where equality holds if and only if $g' = c g$ for some positive
constant $c$}.\\[0.02in]
{\it Proof.} \ \ Let $g = u^{4\over {n - 2}} g'$ with $u > 0$. If $S
(g')$ is the scalar
curvature of the metric $g'$, then
$$C_n \Delta' u - S (g') u = - S (g) u^{{n + 2}\over {n - 2}}\,, \leqno
(4.2)$$
where $C_n = 4 (n - 1)/(n - 2)$ and $\Delta'$ is the Laplacian for the
metric
$g'$. Multiple (4.2) by $u$ and then integrate by parts we
have
\begin{eqnarray*}
- C_n \int_M |\bigtriangledown u|^2_{g'} dv_{g'} - \int_M S (g') u^2
dv_{g'} & = &
|S (g)| \int_M u^{{2n}\over {n - 2}} dv_{g'}\\
& = & |S (g)| {\mbox {Vol}}\, (M, g)\,,
\end{eqnarray*}
as $S (g)$ is a nonpositive constant. Therefore
$$- \int_M S (g') u^2 dv_{g'} \ge |S (g)| {\mbox {Vol}}\, (M, g)\,,
\leqno (4.3)$$
and equality holds if and only if $u$ is a constant.
Using H\"older's inequality we obtain
$$( \int_M |S (g')|^{n \over 2} dv_{g'} )^{2\over n} ( \int_M
u^{{2n}\over {n - 2}} dv_{g'}
 )^{{n - 2}\over n} \ge - \int_M S (g') u^2 dv_{g'}\,.$$
Combine with (4.3) to obtain
$$( \int_M |S (g')|^{n \over 2} dv_{g'} )^{2\over n} ( {\mbox {Vol}}\,
(M, g) )^{{n - 2}\over n}
\ge |S (g)| {\mbox {Vol}}\, (M, g)\,.$$
That is,
$$\int_M |S (g')|^{n \over 2} dv_{g'} \ge |S (g)|^{n\over 2} {\mbox
{Vol}}\, (M, g) =
\int_M |S (g)|^{n \over 2} dv_g\,.$$
{\bf Q.E.D.}\\[0.03in]
\hspace*{0.3in}For a Riemannian metric $g$ on a compact manifold $M$,
the Yamabe invariant is defined as
$$Q (M, g) = \inf \ \{ { { {{4 (n - 1)}\over {n - 2}} \int_M |
\bigtriangledown u |^2 dv_g
+ \int_M R_g u^2 dv_g }\over
{ ( \int_M |u|^{{2n}\over {n - 2}} dv_g )^{{n - 2}\over n} } } \ |
\ \ \ u \in C^\infty (M)\,, \ \ u \not\equiv 0 \ \}\,. \leqno (4.4)$$
It is known that the Yamabe invariant for the standard unit sphere is
equal to the best
constant for the Sobolev inequality on ${\bf R}^n$ (Theorem 3.3 of [14]),
i.e.,
$$Q (S^n, g_o) = n (n - 1) \omega_n^{2\over n}\,,$$
where $\omega_n$ is the volume of the unit $n$-sphere.\\[0.02in]
\hspace*{0.3in}Lemma (4.1) does not hold in general for constant positive
scalar curvature.
However, for Einstein metrics with positive scalar curvature
we have the following result.\\[0.03in]
{\bf Lemma 4.5.} \ \ {\it For $n \ge 3$, let $(M, g_o)$ be a compact
Einstein manifold with positive scalar curvature,
then for any metric $g$ that is conformal to $g_o$, we have}
$$\int_M |S (g)|^{n \over 2} dv_g \ge \int_M |S (g_o)|^{n \over 2}
dv_{g_o}\,.$$

{\it Proof.} \ \ As the scalar curvature of $(M, g_o)$ is positive, we
have $Q (M, g_o) > 0$.
If $Q (M, g_o) < n (n - 1) \omega_n^{2\over n}$, then there is a smooth
positive function $u$ such that
$$Q (M, g_o) = { { {{4 (n - 1)}\over {n - 2}} \int_M | \bigtriangledown u
|^2 dv_g
+ \int_M R_g u^2 dv_g }\over
{ ( \int_M |u|^{{2n}\over {n - 2}} dv_g )^{{n - 2}\over n} } } \,,$$
and the metric $u^{4/(n - 2)} g_o$ has constant positive scalar
curvature. Obata's theorem A implies
that $u$ is a positive constant and
$$Q (M, g_o) = n (n - 1) {\mbox {Vol}}\, (M, g_o)^{2\over n}\,.$$
The same relation holds of the standard $n$-sphere.
(4.4) gives the following
inequality
$$n (n - 1) {\mbox {Vol}}\, (M, g_o)^{2\over n} ( \int_M |u|^{{2n}\over
{n - 2}} dv_{g_o} )^{{n - 2}\over n}
\le 4 {{n - 1}\over {n - 2}} \int_M | \bigtriangledown u |^2 dv_{g_o}
+  \int_M R_{g_o} u^2 dv_{g_o}\,, \leqno (4.6)$$
for $u \in C^\infty (M)\,.$
Let $g = u^{4\over {n - 2}} g_o$, $u > 0$. We have
$$4 {{n - 1}\over {n - 2}} \Delta_o u -  S (g_o) u
= - S (g) u^{{n + 2}\over {n - 2}}\,, \leqno (4.7)$$
where $\Delta_o$ is the Laplacian for $(S^n, g_o)$.
Multiple (4.7) by $u$ and then integrate by parts we
obtain
$$4 {{n - 1}\over {n - 2}} \int_M |\bigtriangledown u|^2 dv_{g_o}
+ \int_M S (g_o) u^2 dv_{g_o} =
   \int_M  S (g) u^{{2n}\over {n - 2}} dv_{g_o}\,. \leqno (4.8)$$
Apply the H\"older's inequality and the inequality (4.6), we have
\begin{eqnarray*}
 \int_M  S (g) u^{{2n}\over {n - 2}} dv_{g_o}
& \le & ( \int_M |S (g)|^{n\over 2} u^{{2n}\over {n - 2}} dv_{g_o}
)^{2\over n}
( \int_M u^{{2n}\over {n - 2}} dv_{g_o} )^{{n - 2}\over n}\\
& \le & [n (n - 1) {\mbox {Vol}}\, (M, g_o)^{2\over n}]^{-1}\,( \int_M |S
(g)|^{n\over 2} dv_{g} )^{2\over n} \,
( 4 {{n - 1}\over {n - 2}} \int_M |\bigtriangledown u|^2 dv_{g_o}  +
  \int_M u^2 dv_{g_o} )\,.
\end{eqnarray*}
So from (4.8) we obtain
\begin{eqnarray*}
& \ & \ \ 4 {{n - 1}\over {n - 2}} \int_M |\bigtriangledown u|^2
dv_{g_o}
+  \int_M S (g) u^2 dv_{g_o}\\
& \ & \le [n (n - 1) {\mbox {Vol}}\, (M, g_o)^{2\over n}]^{-1}\,( \int_M
|S (g)|^{n\over 2} dv_g )^{2\over n}
\,( 4 {{n - 1}\over {n - 2}} \int_M |\bigtriangledown u|^2 dv_{g_o}  +
  \int_M u^2 dv_{g_o} )\,.
\end{eqnarray*}
We must have
$$[n (n - 1) {\mbox {Vol}}\, (M, g_o)^{2\over n}]^{-1}
\,( \int_M |S (g)|^{n\over 2} dv_g )^{2\over n} \ge 1\,,$$
or
$$\int_M |S (g)|^{n\over 2} dv_g \ge [n (n - 1]^{n\over 2} {\mbox
{Vol}}\, (M, g_o)
= \int_M |S (g_o)|^{n\over 2} dv_{g_o}\,,$$
as $S (g_o) = n (n - 1)$.
\ \ \ \ \ {\bf Q.E.D.}\\[0.05in]
{\bf Corollary 4.9.} \ \ {\it  For any metric $g$ on $S^n$
that is conformal to $g_o$ and with $S (g) \le n (n - 1)$, we have {\mbox
{Vol}}\,
$(S^n, g) \ge {\mbox {Vol}}\, (S^n, g_o)$ }\\[0.03in]
{\bf Proposition 4.10.} \ \ {\it  Let $(M, g)$ be an $n$-manifold with $
b^2 g \ge {\mbox {Ric}}\, (g) \ge a^2 g$
for some positive numbers $a$ and $b$. Then for any metric $g' =
u^{4\over {n - 2}} g$, $u > 0$, we have}
$$\int_M |S (g')|^{n \over 2} dv_{g'} \ge c_n
\int_M |S (g)|^{n \over 2} dv_g\,,$$
{\it where $c_n$ is a positive constant that depends on $a$, $b$ and $n$
only.}\\[0.02in]
{\it Proof.} \ \ For the smooth positive function $u$, the Sobolev
inequality on $(M, g)$ [1] gives
$$( \int_M u^{{2n}\over {n - 2}} dv_g)^{{n - 2}\over {2n}}
\le ( {\mbox {Vol}}\, (M, g) )^{-{1\over n}} [ \tau \sigma_n ( \int_M
|\bigtriangledown u|^2 dv_g )^{1\over 2} +
( \int_M u^2 dv_g )^{1\over 2} ]\,, \leqno (4.11)$$
where $\tau = {\mbox {Diam}}\, (M, g) / \alpha_n$ and $\sigma_n$,
$\alpha_n$ are positive constants
that depend on $n$ only. As ${\mbox {Ric}}\, (g) \ge a^2 g\,.$ Myers'
theorem gives Diam $(M, g)
\le \pi {\sqrt {n - 1}} /a$. Therefore there exists a positive constant
$C (n, a)$,
which depends on $n$ and $a$ only,
such that
$$( \int_M u^{{2n}\over {n - 2}} dv_g)^{{n - 2}\over {n}}
\le C (n, a) \,( {\mbox {Vol}}\, (M, g) )^{-{2\over n}} \,( \int_M
|\bigtriangledown u|^2 dv_g  +
 \int_M u^2 dv_g )\,. \leqno (4.12)$$
In the proof of lemma (4.5), if we use the inequality (4.12) instead of
(4.6), we obtain
\begin{eqnarray*}
& \ & 4{{n - 1}\over {n - 2}} \int_M |\bigtriangledown u|^2 dv_g
+ \int_M S (g) u^2 dv_g\\
& \le &
 C (n, a) \,( \int_M |S (g')|^{n\over 2} dv_g )^{2\over n} ( {\mbox
{Vol}}\, (M, g) )^{-{2\over n}} \,( \int_M |\bigtriangledown u|^2 dv_g  +
 \int_M u^2 dv_g )\,.
\end{eqnarray*}
As $S (g) ge n a^2 \ge$, we must have
$$C (n, a) \,( \int_M |S (g')|^{n\over 2} dv_{g'} )^
{2\over n} \,( {\mbox {Vol}}\, (M, g) )^{-{2\over n}} \ge \min\,\,
\{ {{4 (n - 1)}\over {(n - 2)}}, n a^2 \}\,,$$
or
$$\int_M |S (g')|^{n\over 2} dv_{g'}
\ge   \int_M |S (g)|^{n\over 2} dv_g\,,$$
where
$$C (n, a, b) = {{ {\min}\,\,
\{ {{4 (n - 1)}\over {(n - 2)}}, n a^2 \} }\over { C (n, a)
n b^2}}\,.$$
We have made use of the fact that $S (g) \le n b^2$. $C(n, a, b)$ is a
positive constant
that depends on $n$, $a$ and $b$ only. \ \ \ \ \ {\bf Q.E.D.}\\[0.03in]
\hspace*{0.3in}Similar to the Ricci curvature flow, Hamilton has
introduced the normalized
Yamabe flow (scalar curvautre flow):
$${{\partial g_t}\over {\partial t}} = ({\bar s} (g_t) - S (g_t)) g_t\,,
\leqno (4.13)$$
where ${\bar s} (g_t) = \int_M S (g_t) dv_{g_t} / {\mbox {Vol}}\, (M,
g_t)$. The Yamabe flow has
been used by Hamilton, B. Chow [8], and R. Ye [21] to obtain constant
scalar curvature metrics on
various situations. As in section 3, we
consider the change of the $L^{n\over 2}$-norm on scalar curvatures along
the Yamabe flow.\\[0.03in]
{\bf Lemma 4.14.} \ \ {\it Let $(M, g_o)$ be a compact Riemannian
$n$-manifold with $n \ge 4$. Assume that
$(M, g_o)$ has positive
scalar curvature. If $g_t$ is a solution to the Yamabe flow (4.13) with
initial
metric $g_o$, then}
$${d\over {dt}} \int_M |S (g_t)|^{n\over 2} dv_{g_t} \le 0\,,$$
{\it and equality holds at time $t$ if and only if $g_t$ has constant
scalar curvature.}\\[0.02in]
{\it Proof.} \ \ It is more convenient to consider the unnormalized
Yamabe flow
$${{\partial g_t}\over {\partial t}} =  - S (g_t) g_t\,. \leqno (4.15)$$
One can rescale in time for the solutions of (4.15)
to obtain corresponding solutions of (4.13) [8,17].
Under the flow (4.13), the evolution equation for the scalar curvature is [8]
$${{\partial}\over {\partial t}} S (g_t) = (n - 1) \Delta S (g_t) + S
(g_t)^2\,.$$
It follows from the maximal principle that if $g_o$
has positive scalar curvature, then $S (g_t) > 0$ for all $t \ge 0$.
Under the normalized Yamabe flow (4.13),
the evolution equation for the scalar curvature is [18]
$${{\partial}\over {\partial t}} S (g_t)
= (n - 1) \Delta S (g_t) + S (g_t) ( S (g_t) - {\bar s} (g_t) )\,, \leqno
(4.16)$$
and
$$(dv_g)' = {1\over 2} {\mbox {tr}}_g ({{dg}\over {dt}}) dv_g  =
{n\over 2} ({\bar s} (g_t) - S (g_t))\,. \leqno (4.17)$$
Therefore we have
\begin{eqnarray*}
& \ & {d\over {dt}} \int_M |S (g_t)|^{n\over 2} dv_{g_t}\\
& = & \int_M {n\over 2} S (g_t)^{{n\over 2} - 1}
{{\partial}\over {\partial t}} S (g_t) dv_{g_t} + \int_M {n\over 2}
S^{n\over 2} ({\bar s} (g_t) - S (g_t)) dv_{g_t}
 \ \ \ \ ({\mbox {as}} \ \ \ S (g) > 0)\\
& = & \int_M {n\over 2} S (g_t)^{{n\over 2} - 1} [(n - 1) \Delta S (g_t)
+ S (g_t) ( S (g_t) - {\bar s} (g_t) )]
+ \int_M {n\over 2} S^{n\over 2} ({\bar s} (g_t) - S (g_t)) dv_{g_t}\\
& = & - \int_M {n\over 2} ({{n\over 2} - 1}) S (g_t)^{{n\over 2} - 2}
|\bigtriangledown S (g)|^2 dv_{g_t} \le 0\,,
\end{eqnarray*}
and equality holds if and only if $S (g_t)$ is a constant. \ \ \ \ \ {\bf
Q.E.D.}\\[0.05in]
\hspace*{0.3in}Let $(M, g)$ be a compact conformally flat manifold with
positive Ricci curvature. The Yamabe flow (4.6)
with initial metric $g$ is known to converge to a constant curvature
metric $g_o$ as $t \to \infty$ [8]. Applying the above lemma we have
the following.\\[0.03in]
{\bf Theorem 4.18.} \ \ {\it Let $(M, g)$ be a compact conformally flat
manifold with positive Ricci
curvature. Then}
$$\int_M |S (g)|^{n\over 2} dv_{g} \ge \int_M |S (g_o)|^{n\over 2}
dv_{g_o}\,, \leqno (4.19)$$
{\it where $g_o$ has constant positive sectional curvature.}\\[0.05in]
{\it Remark:} \ \ As the Ricci curvature of $(M, g)$ is
positive, it is bounded from below by a positive constant. Hence the
fundemantal
group is finite by Myer's theorem. The universal covering of $M$ is then
conformally
equivalent to the standard $n$-sphere $S^n$ under the development map.
Because a finite group
of conformal trnasformations of the $S^n$ is conjugate to a group of
isometrics of $S^n$, we see that the metric $g$ is conformal to a metric
of $g_o$ of
constant positive sectional curvature. Proposition (4.10) provides a not
so sharp lower
bound on the $L^{n\over 2}$-norm on $S (g)$.\\[0.02in]
\hspace*{0.3in} We note that there exists a family of
metrics on $S^n$ for $n \ge 3$ with $L^{n\over 2}$-norms on the
scalar curvatures concentrate around one point. For
any $\epsilon > 0$, the family of functions
$$u_\epsilon (x) = ( {\epsilon \over {\epsilon^2 + | x |^2}} )^{{n -
2}\over 2}\,, \ \ \  x \in {\bf R}^3\,,$$
satisfy the equation
$$\Delta_o u_\epsilon + n (n - 2) u_\epsilon^{{n + 2}\over {n - 2}} = 0\,,$$
where $\Delta_o$ is the Laplacian for ${\bf R}^n$ with the standard flat
metric
$\delta_{ij}$. That is, the metric $g_{o, \epsilon} = u_\epsilon^{4\over
{n - 2}} \delta_{ij}$ has
scalar curvature equal to $n (n - 2)$. Let $\Phi : S^n \to {\bf R}^n$ be
the sterographic
projection which sends the north pole to infinity. Using the fact that
$d ( (0, 0,...0, 1), y ) \sim 1/ |\Phi (y) |$, where
$(0, 0,.., 0, 1)$ is the north pole of $S^n$ and
$y \in S^n \setminus (0, 0,..., 0, 1)$ and $d$ is the distance on $S^n$,
the pull back of the family of metrics $g_{o, \epsilon}$ by $\Phi$,
denoted by
$g_{\epsilon}$, on
$S^n$, is a family of nonsingular metrics on $S^n$. Then $\Phi :
( S^n \setminus (0, 0,.. 0, 1), g_{\epsilon}) \to ({\bf R}^n\,, g_{o,
\epsilon})$ is
an isometry. The scalar curvature of $(S^n, g_{\epsilon})$ equals to $n
(n - 2)$. And
\begin{eqnarray*}
\int_{S^n} | S (g_{\epsilon}) |^{n\over 2} dv_{g_\epsilon}
& = & \int_{ {\bf R}^n}
[n (n - 2)]^{n\over 2} dv_{g_{o, \epsilon}}\\
& = & \int_{ {\bf R}^n } [n (n - 2)]^{n\over 2} u_{\epsilon}^{{2n}\over
{n - 2}} dv_o\\
& = & \int_{ {\bf R}^n } [n (n - 2)]^{n\over 2}
( {\epsilon \over {\epsilon^2 + | x |^2}} )^n dv_o\\
& = & c_n \int_0^\infty ( {1\over {1 + r^2}} )^n r^{n - 1} dr\,,
\end{eqnarray*}
where $c_n = [n (n - 2)]^{n\over 2} {\mbox {Vol}}\, (S^{n - 1})$ and $r =
|x|/\epsilon, x \in {\bf R}^n$.
As $\epsilon \to 0$, $L^{n\over 2}$-norms on the
scalar curvatures concentrate around the south pole, i.e., there exist a
positive
constant $C_n$ such that
$$\int_{S^n} | S (g_{\epsilon}) |^{n\over 2} dv_{g_\epsilon} \ge C_n$$
for all $1 > \epsilon > 0$ while if $O$ is any open neighborhood of the
south pole, then
$$\int_{S^n \setminus O} | S (g_{\epsilon}) |^{n\over 2} dv_{g_\epsilon}
\to 0 \ \ \
{\mbox {as}} \ \ \ \epsilon \to 0\,.$$
While as $\epsilon \to \infty$, the integral concentrates around the
north pole.

\vspace{0.45in}

{\bf \LARGE Appendix}

\vspace{0.3in}

The two-step variational program purposed to obtain an Einstein metric
on a compact manifold $M$ is to consider the quotient
$$Q (g) = {{ \int_M S (g) \,dv_g}\over { [{\rm {Vol}} (M, g) ]^{{n -
2}\over n}
}}$$
and find a minimizing metric $g_o$ in the conformal class of a metric $g$,
denoted by $C (g)$, that is,
$$\lambda \, (M, g_o) = \inf_{ g' \in C (g)} Q (g') = Q (g_o)\,.$$
($\lambda \, (M, g_o)$ is known as the Yamabe invariant of $M$ for the
conformal class $C (g)$.) Let
$$\Lambda \, (M) = \sup_{g}  \inf_{ g' \in C (g)} Q (g')\,.$$
Then one tries to find a metric
$h$ such that
$$Q (h) = \lambda \, (M, h) = \Lambda \, (M)\,.$$
It is known that if such a Riemannian metric $h$ exists, then it is an
Einstein
metric. Let $M$ be a compact manifold which admits a metric of negative
sectional curvature. Then by the solution to the Yamabe problem, on each
conformal class $C (g)$ of Riemannian metrics on $M$, there is a unique
metric
$g_o \in C (g)$ with $S (g_o) = - 1$ and
$$\lambda \, (M, g_o) = - [{\rm {Vol}} (M, g_o) ]^{{2}\over n}\,.$$
Hence in this situation $\Lambda \, (M)$ is related to a lower bound on
Vol $(M,
g)$ for any Riemannian metric with $S (g) = -1$. This is equivalent to an
lower
bound on the $L^{n/2}$-norm of scalar curvature, as follows from [3] that
$$\int_g |S (g)|^{n\over 2} dv_g \ge \int_{g_o} |S (g_o)|^{n\over 2}
dv_{g_o}$$
for any $g \in C (g_o)$ and $S (g_o) \equiv -1$.\\[0.3in]

{\bf Proposition.} {\bf  Let $M$ be a compact complex surface which
admits a
K\"ahler-Einstein metric $\eta$ of negative scalar curvature. Then up to
scaling and isometry,
$\eta$ is the unique metric on $M$ such that}
$$Q (\eta) = \lambda \, (M, \eta) = \Lambda \, (M)\,.$$

{\it Proof.} \ \  From theorem 2.12 and the above remark, it follows that
$$Q (\eta) = \lambda \, (M, \eta) = \Lambda \, (M)\,.$$
Suppose that $g$ is another metric on $M$ such that
$$Q (g) = \lambda \, (M, g) = \Lambda \, (M)\,.$$
Then $g$ is an Einstein metric with negative scalar curvature, as
$\Lambda \,
(M)$ is negative [16]. We may assume that, without loss of generality, $S
(\eta) = S
(g) = -1$. Then
$$Q (\eta) = \Lambda \, (M) = Q (g)$$
implies that
$$\int_M |S (g)|^2 dv_{g} = \int_M |S (\eta)|^2 dv_\eta$$
and
$$c_1 (L)^2 = {1\over {32 \pi^2}} \int_M |S (g)|^2 dv_g\,.$$
A result in [13] implies that $g$ is a K\"ahler-Einstein metric. Hence
$g$ is
isometric to $\eta$ up to a scaling factor [3]. \ \ \ \ \ {\bf
Q.E.D.}\\[0.03in]

\hspace*{0.3in}For a compact hyperbolic $4$-manifold $(M, h)$, the
hyperbolic metric
$h$ is the only candidate to achieve
$$Q (h) = \lambda \, (M, h) = \Lambda \, (M)\,.$$
For if there is another Riemannian metric $g$ such that
$$Q (g) = \lambda \, (M, g) = \Lambda \, (M)\,,$$
then $g$ is an Einstein metric. Hence $g$ is isometric to $h$ up to
scaling, by
a result of Besson-Courtois-Besson [4, 5]. However, it is unsettled
whether
$Q (h) = \Lambda \, (M)$ (see [16]).

\pagebreak

\vspace{0.3in}

\centerline{\bf \LARGE References}

\vspace{0.3in}

\ [1] P. Berard, {\it From vanishing theorems to estimating theorems: the
Bochner technique revisited,}
Bull. Amer. Math. Soc. {\bf 19} (1988), 371-406.\\

\ [2] M. Berger, P. Gauduchon \& E. Mazet, {\it Le Spectre d'Une
Vari\'et\'e
Riemannienne,} Lecture Notes in Math. {\bf 194}, Springer-Verlag,
Berlin-New York, 1971.\\

\ [3] A. Besse, {\it Einstein Manifolds,} Springer-Verlag, Berlin-New
York,
1987.\\

\ [4] G. Besson, G. Courtois \& S. Gallot, {\it Volume et entropie
minimale des espaces localement sym\'etriques,}
Invent. Math. {\bf 103} (1991), 417-445.\\

\ [5] G. Besson, G. Courtois \& S. Gallot, {\it Les vari\'et\'es
hyperboliques sont des minima locaux de l'entropie topologique,}
Invent. Math. {\bf 117} (1994), 403-445.\\

\ [6] G. Besson, G. Courtois \& S. Gallot, {\it Entropies et Rigidit\'es des
Espaces Localement Sym\'etriques de Courbure Strictement N\'egative,}
Pr\'epublication de l'Institut Fourier, Grenoble, No. 281, (1994).\\

\ [7] J. Cheeger \& D. Ebin, {\it Comparison theorems in Riemannian
geometry,}
North-Holland Mathematical Library, vol. 9. \\

\ [8] B. Chow, {\it The Yamabe flow on locally conformally flat manifolds
with positive Ricci curvature,}
Comm. Pure and Applied Math. {\bf 45} (1992), 1003-1014.\\

\ [9] L. Z. Gao, {\it Convergence of Riemannian manifolds; Ricci and
$L^{n/2}$-curvature pinching,}
J. Diff. Geom. {\bf 32} (1990) 349-381.\\

[10] M. J. Gursky, {\it Locally conformally flat four- and six-manifolds
of positive
scalar curvature and positive Euler characteristic,} Indiana Univ. Math.
J. {\bf 43}
(1994), 747-774.\\

[11] R. Hamilton, {\it Three-manifolds with positive Ricci curvature,} J.
Diff. Geom. {\bf 17}
(1982), 225-306.\\

[12] E. Hebey \& M. Vaugon, {\it Un th\'eor\'eme de pincement int\'egral
sur la courbure
concirculaire en g\'eom\'etrie conforme,} C.R. Acad. Sci. Serie I {\bf
316} (1993), 483-488.\\

[13] C. LeBrun, {\it Einstein metrics and Mostow rigidity,} Math.
Research Letters
{\bf 2} (1995), 1--8.\\

[14] J. M. Lee \& T. H. Parker, {\it The Yamabe problem,} Bull. Amer.
Math. Soc. {\bf 17} (1987), 37-91.\\

[15] M. Min-Oo, {\it Almost Einstein manifolds of negative Ricci
curvature,} J. Diff. Geom.
{\bf 32} (1990), 457-472.\\

[16] R. Schoen, {\it Variational theory for the total scalar curvature
functional
for Riemannian metrics and related topics,} Lecture Notes in Math. {\bf
1365}, 120-154.\\

[17] Z. Shen, {\it Some rigidity phenomena for Einstein metrics,} Proc.
Amer. Math. Soc.
{\bf 108} (1990), 981-987.\\

[18] M. Spivak, {\it A comprehensive introduction to differential
geometry,} vol. 5,
Publish or Perish, 1975.\\

[19] E. Witten, {\it  Monopoles and four-manifolds,} Math. Research Letters
{\bf 1} (1994), 769--796.\\

[20] R. Ye, {\it Ricci Flow, Einstein metrics and space forms,} Trans.
Amer. Math. Soc.
{\bf 338} (1993), 871-896.\\

[21] R. Ye, {\it Global existence and convergence of Yamabe flow,} J.
Diff. Geom. {\bf 17}
(1994), 35-50.\\

\end{document}